\documentclass[12pt]{article}
\usepackage[centertags]{amsmath}
\usepackage{amsfonts}
\usepackage{amssymb}
\usepackage{amsthm}
\usepackage{amsmath,mathrsfs}
\usepackage{amsmath,amscd}
\usepackage{mathtools}
\usepackage[matrix,arrow,curve,cmtip]{xy}
\title{\textbf{The Time-Energy Uncertainty Principle in Algebraic Geometry}}
\author{Renaud Gauthier \footnote{2020 Math. Subj. Class: 14A30, 81P17, 81S07 . Keywords: derived stacks, time-energy uncertainty relation, derivations, energy, Shannon's entropy} \\ \\}
\theoremstyle{definition}

\newtheorem{EAK}{Lemma}[subsubsection]
\newtheorem{EpresHcov}[EAK]{Proposition}
\newtheorem{EpresHFP}{Proposition}[subsubsection]

\newcommand{\beq}{\begin{equation}}
\newcommand{\eeq}{\end{equation}}

\newcommand{\rarr}{\rightarrow}

\newcommand{\larr}{\leftarrow}

\newcommand{\xrarr}{\xrightarrow}

\newcommand{\cA}{\mathcal{A}}

\newcommand{\cC}{\mathcal{C}}

\newcommand{\cD}{\mathcal{D}}
\newcommand{\cE}{\mathcal{E}}

\newcommand{\cG}{\mathcal{G}}
\newcommand{\cH}{\mathcal{H}}

\newcommand{\cK}{\mathcal{K}}

\newcommand{\cO}{\mathcal{O}}

\newcommand{\cX}{\mathcal{X}}

\newcommand{\bR}{\mathbb{R}}

\newcommand{\Cat}{\text{Cat}}

\newcommand{\Fun}{\text{Fun}}

\newcommand{\Hom}{\text{Hom}}
\newcommand{\Ho}{\text{Ho}\,}

\newcommand{\Map}{\text{Map}}
\newcommand{\map}{\text{map}}

\newcommand{\op}{\text{op}}

\newcommand{\Set}{\text{Set}}
\newcommand{\Spec}{\text{Spec\,}}
\newcommand{\Top}{\text{Top}}
\newcommand{\uHom}{\underline{\Hom}}

\newcommand{\AMod}{A\text{-Mod}}

\newcommand{\AffC}{\cA \text{ff}_{\cC}}
\newcommand{\AffChat}{\AffC^{\,\wedge}}
\newcommand{\AffCtildetau}{\AffC^{\,\sim,\tau}}

\newcommand{\cHom}{\cH \text{om}}

\newcommand{\Comm}{\text{Comm}}

\newcommand{\Catinf}{\Cat_{\infty}}

\newcommand{\CommC}{\text{Comm}(\cC)}

\newcommand{\Dn}{\Delta^n}

\newcommand{\Der}{\mathbb{D}\text{er}}

\newcommand{\dSt}{\text{dSt}}

\newcommand{\eps}{\epsilon}

\newcommand{\oP}{\oplus}

\newcommand{\oT}{\otimes}

\newcommand{\pDn}{\partial \Delta^n}

\newcommand{\RpDn}{\mathbb{R}\pDn}
\newcommand{\RDn}{\mathbb{R}\Dn}

\newcommand{\RuHom}{\bR \uHom}
\newcommand{\RuSpec}{\bR \uSpec}

\newcommand{\sPr}{\text{sPr}}

\newcommand{\SetD}{\Set_{\Delta}}

\newcommand{\skMod}{\text{s}k\text{-Mod}}

\newcommand{\skAlg}{\text{s}k\text{-Alg}}

\newcommand{\uh}{\underline{h}}

\newcommand{\uK}{\underline{K}}

\newcommand{\ucHom}{\underline{\cHom}}

\newcommand{\uSpec}{\underline{\Spec}}

\newcommand{\AffCE}{\cA \text{ff}_{\cC, E}}
\newcommand{\tauE}{\tau_E}
\newcommand{\AffCEtildetauE}{\AffCE^{\,\sim,\tauE}}
\newcommand{\CE}{\cC_E}
\newcommand{\Cov}{\text{Cov}}
\newcommand{\dAffC}{\text{d}\AffC}
\newcommand{\dAffChat}{\dAffC^{\wedge}}
\newcommand{\dAffCtildetau}{\text{d}\AffCtildetau}

\begin{document}
\maketitle
\begin{abstract}
We consider the time energy uncertainty principle from Quantum Mechanics and provide its Algebro-Geometric interpretation within the context of stacks.
\end{abstract}

\section{Introduction}
In the papers \cite{RG}, \cite{RG2}, \cite{RG3}, \cite{RG4} we have modeled natural phenomena using  derived stacks. To be precise, we worked with $\cX = \dSt(k)$ a Segal category of derived stacks, with a base symmetric monoidal model category $\cC = \skMod$ of simplicial modules over a commutative ring $k$, $\CommC = \skAlg$ and $\dAffC = L\CommC^{\op}$, a Dwyer-Kan simplicial localization, $\dSt(k) = \dAffCtildetau$ being a left Bousfield localization of $\dAffChat = \RuHom(\dAffC^{\op}, \Top)$, with $\Top = L \SetD$, $\tau$ a Segal topology on $\dAffC$, or equivalently $\cX$ can be obtained as a simplicial localization of $\sPr_{\tau}(\dAffC) \simeq \AffCtildetau$, the model category of stacks on $\AffC$, itself a double Bousfield localization of $\sPr(\AffC)$, the category of simplicial presheaves on $\AffC$. In spite of the fact that $\cX$ is really the Segal category of stacks we had in mind in the above references, presently we consider a generic symmetric monoidal model category $\cC$ with simplicial objects (something we need for the time-energy uncertainty relation), satisfying the properties of a homotopical algebraic context (\cite{HAGII} - to develop our theory within the context of HAG), and we work with the corresponding model category of stacks $\AffCtildetau$ for the sake of generality, $\tau$ being a model topology on $\AffC$. We will focus on those model categories of stacks that model natural phenomena. We argue functors $F$ in such model categories can be decomposed as a composite of functors $F = V \circ E$, $E$ being a stack in $\AffCtildetau$ playing the role of energy, and from that perspective has for sole purpose to materialize the model category $\AffC$ in the natural realm, thereby producing a replica $\AffCE$ of $\AffC$, a natural enrichment. Note that this means $E$ preserves equivalences and small limits, in particular pullbacks, something we will use. Functorially we initially ask that $E$ preserves covers, as a consequence of which it will also preserve homotopy fiber products, something we refer to as being homotopy continuous. For the sake of a given site $(\AffCE, \tauE)$ being an enrichment of the site $(\AffC,\tau)$, we also need that $E$ reflects covers (in addition to preserving them) in the sense that $\{ EX \rarr EY \} \in \Cov(EY)$ in $\Ho(\AffCE)$ implies $\{X \rarr Y\} \in \Cov(Y)$ in $\Ho(\AffC)$. Finally we ask that $E$ preserves functorial factorizations and is $\SetD$-equivariant. It will follow $E$ also preserves hypercovers. From there one considers objects $V \in \AffCEtildetauE$, which we refer to as observables, and which describe dynamics, hence can be used to define a notion of time. Such a factorization $V \circ E$ we call a brane, characterized by an energy $E$, and a pulsation $V$ that generates dynamics. More precisely we define a category of $E$-branes to be a pair $(E,\AffCEtildetauE)$ with objects pairs $(E,V)$, where $V \in\AffCEtildetauE$ are the observables just mentioned. It is when one considers the time-energy uncertainty principle that we will ask that $E$ be a stack, so at that point the only additional requirements on $E$ we import from the above discussion that are not covered by the stack condition are preserving and reflecting covers, preserving functorial factorizations, and being $\SetD$-equivariant.\\

Regarding the time-energy uncertainty relation, one aims for a relation of the form $\Delta E \Delta t \geq \hbar/2$. We will see that $t$ as a variable is not well-defined, but is implicitly defined through the variations of observables $V$, so one needs to make sense of the variations $\Delta E$ and $\Delta V$. We argue the functor of derivations in the study of the infinitesimal theory of stacks provides a functorial way to deal with variations of this nature. To wit, if $A \in \CommC$, $M \in \AMod$, $F$ is a stack, then $F(A \oplus M) \rarr F(A)$ provides a deformation of $F(A)$ in the $M$-direction. The homotopy fiber of this map provides a lift of the simplicial set of derivations $\Der_F(X,M)$, where $X = \RuSpec A$, $\RuSpec = \bR \uh: \Ho(\AffC) \rarr \Ho(\AffChat)$ is the derived model Yoneda. Doing this for all $M$ generates $\Delta F(A)$. Functorially this is given by $\Der_F(X,-)$. Doing this for any $X$ gives us a bifunctor $\Delta F = \Der_F$. Applying this to the functors $E$ and $V$ one can make sense of $\Delta E$ and $\Delta V$ individually. According to \cite{M}, the time-energy uncertainty relation reads:
\beq
\frac{\Delta V}{|\frac{d\langle V \rangle}{dt}|} \Delta H \geq \hbar/2 \nonumber
\eeq
where for us $E$ will play the role of the hamiltonian $H$, and $V$ above is an observable. We argue the quotient $\Delta V/|\frac{d\langle V \rangle}{dt}|$ indicates that one should express $\Delta V = \Der_V$ relative to the tangent stack $TV$. There is a section $V \rarr TV$ that induces a morphism $\Der_V(Y,-) \rarr \Der_{TV}(Y,-)$ for $Y = \RuSpec B$, $B \in \Comm(\cC_E)$, whose standard homotopy fibers generate the functor of relative derivations $\Der_{V/TV}(Y,-)$, and we consider the associated bifunctor $\Der_{V/TV}$, which represents the desired fraction. In a non-quantum algebra treatment of the time-energy uncertainty relation, we therefore ask that the product functor $\Der_{V/TV} \times \Der_E$ be non-contractible, and this provides the Algebro-geometric version of the time-energy uncertainty relation.\\

\section*{Notations}
We will use \cite{Lu} and \cite{Lu2} extensively, so for the sake of simplicity we will refer to those references as HTT and HA respectively.

\section{Energy}
\subsection{Energy as realization functor}
Information, such as that encapsulated by algebraic laws, is pure intellection and is not manifestly realized in nature, hence for the sake of studying algebra induced natural phenomena one sees the need to introduce a realization map to manifest information of this type in nature, and we designate by the name energy, for obvious reasons, that necessary map. Energy operates the transfer from pure information to physical event. An easy way to implement this is to see energy as a functor from a category that we regard as being purely information theoretic in nature, into another one that is a reflexion thereof. More precisely if $E$ denotes energy, we regard it as a functor $E: \cC \rarr \CE$, where if $\cC$ is an abstract category, $\CE$ is its realization in nature, an enhanced replica of $\cC$. In practice we will ask that $\cC$ be a symmetric monoidal category, hence so is $\CE$, since it is an enhanced version of $\cC$, so that $E$ itself is a symmetric monoidal functor.\\

To fix notations, as argued in the introduction, we work in full generality, and for that reason even though one has $E \in \dSt(k)$ in practice, one considers $E$ as a stack in $\AffCtildetau$, for $\cC$ a symmetric monoidal model category with simplicial objects satisfying the conditions of a homotopical algebraic context, and $\tau$ is a model topology on $\AffC$, following the notations of \cite{HAGII}.\\

\subsection{Properties of energy}
We have seen that derived stacks $F \in \dSt(k)$, or more generally objects of $\AffCtildetau$, describe natural phenomena, so they ought to materialize algebraic laws in the natural realm. Energy actually does just that we postulate, hence stacks must factor through energy. Now among algebraic laws in $\AffC$, some govern natural phenomena, or equivalently said, dynamics, thus the stacks that make those laws manifest in nature are first factored through energy, followed by a functor that describes dynamics, and observables $V$ take care of that, leading to a factorization $F = V \circ E$ of stacks in general if one allows that $V = id$ exists. Presently we mention the observables in such decompositions since many of the properties satisfied by energy will follow from having such a factorization of stacks $F \in \AffCtildetau$ in general.\\

\subsubsection{Minimal preservations}
For the sake of interpreting the time-energy uncertainty principle in Algebraic Geometry, we argue $V$ in the factorization $V \circ E$ is a stack, and for that to be true $E$ itself needs to produce an intermediate site $(\AffCE,\tauE)$ on which observables $V$ are defined. We ask that an energy $E$ merely reflects $(\AffC,\tau)$ into a site $(\AffCE,\tauE)$, since for one thing energy materializes algebraic laws in nature, hence $\AffCE$ is a replica of $\AffC$ (which means $E$ preserves equivalences and small limits such as pullbacks), but dynamically one also wants to preserve topologies, so $E(\tau) = \tauE$, and in particular that means \underline{energy preserves and reflects covers}. Lastly, later we will ask that $E$ preserves hypercovers, and for this to be true we will require \newline \underline{$E$ to preserve functorial factorizations and to be $\SetD$-equivariant}. What this means is that if $X \rarr \cE(X) \rarr Y$ is a functorial factorization in $\AffC$, this maps under $E$ to $EX \rarr E(\cE(X)) \rarr EY$, and to ask that $E$ preserves functorial factorizations means this is equivalent to $EX \rarr \cE_E(EX) \rarr EY$, $\cE_E$ functorial factorization in $\AffCE$, in the sense that for all $X \in \AffC$, $E(\cE X) \xrarr{\sim} \cE_E(E(X))$ in $\AffCE$, which can easily be implemented by asking that the fixed functorial factorization on $\AffCE$ is dictated by the one we have on $\AffC$, that is $\cE_E(E) = E(\cE)$.

\subsubsection{Topology $\tauE$ on $\AffCE$}
$E(\cC) = \cC_E$ being an $E$-copy of $\cC$ in the enriched sense, we extend this characteristic of $E$ to the topology on $\AffC$, and we ask that $\tauE$ on $\AffCE$ be a replica of $\tau$ on $\AffC$. A minimal requirement for this to be true is that $E$ preserves covers, and for covers in $\tauE$ to exactly come from those of $\tau$ on $\AffC$ if we are to talk about an enrichment, we ask that $E$ reflects covers, in the sense that if $ EX \rarr  EY$ is a covering map in $\Ho(\AffCE)$, then $X \rarr Y$ is a covering map in $\Ho(\AffC)$. One could have simply asked that $\tauE$ be an $E$-enrichment of $\tau$. This is equivalent to saying that $E$ preserves and reflects covers, and as we will see it must preserve pullbacks as well, which it does, and those are properties we will use later in another context, so we focus on those conditions of $E$.\\

To justify that $\tauE$ is well-defined, as a Grothendieck topology on $Ho(\AffCE)$, if $EX \rarr EY$ is an isomorphism in $\Ho(\AffCE)$, $\cC_E$ being a copy of $\cC$, it must come from an isomorphism in $\Ho(\AffC)$, so is a covering map there, and $E$ preserves covers, so $EX \rarr EY$ is a covering map as well. If $\{EX_i \rarr EX\} \in \Cov(EX)$, that is $\{X_i \rarr X \} \in \Cov(X)$ in $\Ho(\AffC)$, if $\{EX_{ij} \rarr EX_i\} \in \Cov(EX_i)$ in $\Ho(\AffCE)$, that is $\{X_{ij} \rarr X_i\} \in \Cov(X_i)$ in $\Ho(\AffC)$, then $\{X_{ij} \rarr X \} \in \Cov(X)$ in $\Ho(\AffC)$, and $E$ preserves covers, so $\{EX_{ij} \rarr EX\} \in \Cov(EX)$ in $\Ho(\AffCE)$, thus we have stability under composition. For pullbacks, if $\{EX_i \rarr EX \} \in \Cov(EX)$ and $EY \rarr EX$ is any map, we want $\{EX_i \times_{EX} EY \rarr EY \} \in \Cov(EY)$. Since $E$ preserves pullbacks (this is all we need, preserving small limits in full is not necessary) so $EX_i \times_{EX} EY = E(X_i \times_X Y) \rarr EY$, thus $\{E(X_i \times_X Y) \rarr EY \} \in \Cov(EY)$, so that $\{X_i \times_X Y \rarr Y \} \in \Cov(Y)$, and $E$ preserves covers, so this maps back to $\{E(X_i \times_X Y) = EX_i \times_{EX} EY \rarr EY \} \in \Cov(EY)$ in $\Ho(\AffCE)$.

\subsubsection{Energy preserves homotopy fiber products}
Additionally for any finite family $\{A_i\}$ in $\Comm(\cC)$, using the fact that $F$ and $V$ are stacks:
\beq
\xymatrix{
	F(\prod^h A_i) = V \circ E(\prod^h A_i) \ar@{.>}[rr]^{\cong} \ar[dr]_{\cong}
	&& V(\prod^h EA_i)  \ar[dl]_{\cong} \\
	 &\prod FA_i = \prod VEA_i
} \nonumber
\eeq
The horizontal isomorphism would follow as a consequence of $E$ preserving homotopy fiber products. We show this holds:
\begin{EpresHFP}
	Energy preserves homotopy fiber products
\end{EpresHFP}
\begin{proof}
	For any map $X \rarr Y$ in $\AffC$, fix a functorial factorization $X \rarr \cE(X) \rarr Y$, in such a manner that $\cE(X) \times_Z \cE(Y) = X \times^h_Z Y$ for maps $X \rarr Z$ and $Y \rarr Z$. We have already assumed that energy $E$ preserves functorial factorizations, which means that $E(\cE) = \cE_E(E)$, and it preserves pullbacks. It follows:
	\begin{align}
		E(X \times^h_Z Y) &= E \big (\cE(X) \times_Z \cE(Y) \big ) \nonumber \\
		&=E(\cE(X)) \times_{EZ} E(\cE(Y)) \nonumber \\
		&=\cE_E(EX) \times_{EZ} \cE_E(EY) \nonumber \\
		&=EX \times^h_{EZ} EY \nonumber
	\end{align}
\end{proof}

\subsubsection{Energy preserves hypercovers}
Additionally for any hypercover $X_* = \Spec A_* \rarr Y = \Spec B$ in $(\AffC,\tau)$, with associated co-simplicial augmented map $A_* \larr B$, $F$ being a stack we know $FB \rarr \lim FA_n$ is an equivalence. Since $V$ is a stack in the decomposition $F = V \circ E$, and $F(B) = V(EB) \rarr \lim V(EA_n) = \lim F(A_n)$ being an equivalence, this suggests regarding $ \Spec EA_* \rarr \Spec EB$ as a hypercover in $(\AffCE,\tauE)$ with associated map of algebras $EA_* \larr EB$, which prompts us to show $E$ preserves the co-simplicial augmented maps associated with hypercovers, which is equivalent to saying $E$ preserves hypercovers. We make this precise.\\

Following \cite{HAGII}, regarding covers, if $(M,\tau)$ is a site, a morphism $X_* \rarr Y_*$ in $sM$ is a $\tau$-hypercover if for all n:
\beq
X_*^{\RDn} \simeq X_n \rarr X_*^{\RpDn} \times^h_{Y_*^{\RpDn}} Y_*^{\RDn} \nonumber
\eeq
is a covering in $\Ho(M)$. Here $X_*^K = (X_*^{\uK})_0$, where the exponential object $X_*^{\uK}$ is defined by the adjunction formula $\Hom(\uK \oT Y_*,X_*) \cong \Hom(Y_*,X_*^{\uK})$ and we consider that simplicial objects of $M$ are tensored and co-tensored in $\SetD$, by defining $\uK \oT Y_*$ as the functor $[n] \mapsto \coprod_{K_n} Y_n$, where $K \in \SetD$. We apply this to the case $M = \AffC = \CommC^{\op}$.\\

One of the things that we will ask of $E$ is that it be \textbf{$\SetD$-equivariant}, which means that if $A_* \in s\CommC$:
\beq
E(\uK \oT A_*) \cong \uK \oT EA_* \nonumber
\eeq
Observe also that if $A_*$ is a simplicial object of $\Comm(\cC)$, $EA_* \in s\Comm(\CE)$. In particular:
\beq
(E(A_*^{\uK}))_0 \cong E((A_*^{\uK})_0) \nonumber
\eeq

\begin{EAK} 
For $A_* \in s\CommC$, we have	$E(A_*^K) \cong (EA_*)^K$.
\end{EAK}
\begin{proof}
We first show $E(A_*^{\uK}) \cong E(A_*)^{\uK}$. We start with:
\beq
\Hom(\uK \oT B_*,A_*) \cong \Hom(B_*,A_*^{\uK}) \nonumber
\eeq
which is mapped by $E$ to:
\beq
	\Hom(E(\uK \oT B_*),EA_*) \cong \Hom(EB_*,E(A_*^{\uK})) \nonumber
\eeq
	However as we argued $E$ is $\SetD$-equivariant so $E(\uK \oT B_*) \cong \uK \oT EB_*$ so the LHS produces:
	\begin{align}
		\Hom (\uK \oT EB_*,EA_*) \cong \Hom(EB_*, (EA_*)^{\uK}) \nonumber
	\end{align}
From which we get $E(A_*^{\uK}) \cong EA_*^{\uK}$. Taking the zero-th component on both sides, we get:
	\begin{align}
		(EA_*)^K = ((EA_*)^{\uK})_0 &\cong (E(A_*^{\uK}))_0 \nonumber \\
		& \cong E((A_*^{\uK})_0)  \nonumber \\
		& = E(A_*^K) \nonumber
	\end{align}
\end{proof}

\begin{EpresHcov}
	Energy $E$ preserves hypercovers
\end{EpresHcov}
\begin{proof}
Consider a hypercover $X_* = \Spec A_* \rarr Y = \Spec B$ in $\AffC$, that is for all $n$, $X_*^{\RDn} \simeq X_n \rarr X_*^{\RpDn} \times^h_{Y^{\RpDn}} Y$ is a covering in $\Ho(\AffC)$. $E$ preserves coverings, and homotopy fiber products, so this maps to:
\begin{align}
	EX_*^{\RDn} \simeq EX_n &\rarr E(X_*^{\RpDn} \times^h_{Y^{\RpDn}} Y) \nonumber \\
	&\cong E(X_*^{\RpDn}) \times^h_{E(Y^{\RpDn})} EY \nonumber
\end{align}
	which is again a covering map. But we have shown $E(A_*^K) \cong (EA_*)^K$, thus with $X_* = \Spec A_* = A_*^{\op} = A_*$ we have:
	\beq
	E(X_*^K) = E(A_*^K) \cong (EA_*)^K = (EX_*)^K \nonumber
	\eeq
	from which it follows:
\beq
(EX)_n \rarr (EX_*)^{\RpDn} \times^h_{(EY)^{\RpDn}} EY \nonumber
\eeq
is a covering in $\Ho(\AffCE)$ for all $n$, hence $EX_* \rarr EY$ is a hypercover, as was to be shown.
\end{proof}

Observe that what we have shown is that if $X_* = \Spec A_* \rarr Y = \Spec B$ is a hypercover, so is $EX_* = E\Spec A_* \rarr EY = E\Spec B$, while for $V$ being a stack we need $V(EB) \xrarr{\simeq} \lim V(EA_n)$, which would follow from $\Spec EA_* \rarr \Spec EB$ being an hypercover in $(\AffCE, \tauE)$. We show those are the same.
\beq
\Spec EA_* = (EA_*)^{\op} = EA_* = EA_*^{\op} = E\Spec A_* = EX_* \rarr EY = E \Spec B \nonumber
\eeq
To express this in terms of algebras, starting then from hypercovers in $\AffC$, a hypercover $X_* = \Spec A_* \rarr Y = \Spec B$ corresponds to a co-simplicial object $A_* \larr B$ in $\Comm(\cC)$, mapping to $EA_* \larr EB$ in $\Comm(\CE)$, which corresponds to $\Spec EA_* \rarr \Spec EB$ in $\AffCE$, which is again a hypercover as just argued. At this point the descent condition for $V \in \AffCEtildetauE$ reads $V(EB) \rarr \lim V(EA_n)$ is an equivalence, which one can rewrite in the form $F(B) \xrarr{\sim} \lim F(A_n)$, where $F = V \circ E$.\\

\subsubsection{Energy as a stack}
To summarize, energies $E$ preserve covers and homotopy fiber products (as a direct consequence of preserving functorial factorizations and products), something we refer to as being \textbf{homotopy continuous}, and they are also $\SetD$-equivariant and reflect covers. These are minimal requirements we put on energies $E$. It follows from these assumptions that energies $E$ preserve topologies, equivalences, homotopy fiber products as just seen and hypercovers. Those properties must hold for $E$ for modeling stacks to factor through energy. However for the sake of defining variations $\Delta E$ involved in the time-energy uncertainty principle in an algebro-geometric setting, we also need that $E$ be a stack, so preserving equivalences, homotopy fiber products and hypercovers is just part of that definition, and the only additional things we ask of energy is that they preserve and reflect covers, functorial factorizations and be $\SetD$-equivariant. Observe that to ask that the objects of $\cC$ be simplicial justifies that $E$ is valued in $\SetD$, which is needed since stacks are localizations of simplicial presheaves.\\

\section{Uncertainty}
We designate by $\Delta E$ the amount of focus put on $E$. To be precise, for $A \in \CommC$, a perfect measurement of energy amounts to knowing $E(A)$ precisely. For possible deviations $A \oplus M$ away from $A$ (contributing to enlarging our neighborhood of $A$), where $M \in \AMod$, we have corresponding objects $E(A \oplus M)$, which consequently contribute to a lack of focus, so from that perspective one sees objects of this nature to contribute to $\Delta E(A)$. We will argue the collection of maps of the form $E(A \oP M) \rarr E(A)$ defines $\Delta E(A)$. A proper way to define $\Delta E$ algebraically goes as follows. Recall from \cite{HAGII} that if $F \in \AffCtildetau$, $A \in \CommC$, $x: X = \RuSpec A \rarr F$, $\uSpec A = \uh_{\Spec A}$, $\uh: \AffC \rarr \AffChat$ the model Yoneda embedding, $M \in \AMod$, $X[M] = \RuSpec(A \oplus M)$, $\Der_F(X,M) = \Map_{X/\AffCtildetau}(X[M],F) \in \Ho(\SetD)$ is the simplicial set of derivations from $F$ to $M$ at the point $x$, which naively reads $\Map(\RuSpec(A \oP M),F) \simeq F(A \oP M)$. Then observe that the functor:
\begin{align}
	\cD: \AMod & \rarr \SetD \nonumber \\
	M &\mapsto hofiber(F(A \oplus M) \rarr F(A)) = \cD(M) \nonumber
\end{align}
appears as a lift of $\Der_F(X,-): \Ho(\AMod) \rarr \Ho(\SetD)$. However if $\cD$ represents a satisfying representation of $\Delta E(A)$ at the level of model categories, we regard $\Der_E(X,-)$ as the appropriate functorial representation of $\Delta E(A)$ at the level of categories, and $\Der_E(X,-)$ being functorial in $X$, $\Delta E \doteq \Der_E$ is the algebro-geometric representation of $\Delta E$.

\section{Branes}
Given $E$, we have observables $V \in \AffCEtildetauE$ that correspond to pulsations of natural phenomena generated by $E$. We regard the pair $(E,\AffCEtildetauE)$ as a category of branes, which we will refer to as $E$-branes, with objects $(E,V)$ that simply provide a factorization of objects $F$ of $\AffCtildetau$ as $F = V \circ E$, with $E \in \AffCtildetau$ materializing algebraic laws in the world of natural phenomena, and $V$ providing the dynamics, hence a notion of time. In what follows branes will be used to illustrate the time-energy uncertainty relation.\\

\section{Time-energy uncertainty principle}
Recall from \cite{M} that if $A$ is an observable and $\psi$ is a state of our system, we write $\langle A \rangle = \langle \psi | A | \psi \rangle$. We have $\Delta A = \sqrt{\langle \psi | (A - \langle A \rangle )^2 | \psi \rangle}$, thus $\Delta A$ appears as a deformation of $A$ away from $\langle A \rangle$ in the state $\psi$. Now for $H$ the hamiltonian, we have $\Delta A \Delta H \geq \hbar/2 | \frac{d \langle A \rangle}{dt}|$, which one can rewrite in the form:
\beq
\frac{\Delta A}{|\frac{d\langle A \rangle}{dt}|} \Delta H \geq \hbar/2 \nonumber
\eeq
something that we can naively rewrite as $dt \Delta H \geq \hbar/2$, hence the old fashioned time-energy uncertainty relation. What this means is that $\Delta A$ is measured in units of $|\frac{d\langle A \rangle}{dt}|$. Additionally, we see that time itself is not a standalone variable, but arises from dynamics; one considers the focus put on a certain observable $A$, and $\Delta A$ is expressed in terms of a tangent object, resulting in a notion of time differentials. This being true for a predetermined $H$, in an algebro-geometric setting $E \in \AffCtildetau$ plays the role of $H$, and observables $V \in \AffCEtildetauE$ represent $A$ above.\\

Presently we focus on $|\frac{d\langle A \rangle}{dt}|$. If $A$ is an observable, represented as a stack then, $\langle A \rangle$ is functorially represented as $A(\psi)$ (which in our situation reads $V(B)$ for $B \in \Comm(\cC_E)$). Computing the derivative of such an object amounts to working with $TA$, its associated tangent stack. Recall that the tangent stack $TF$ of a stack $F$ is defined as $\mathbb{R}_{\tau} \ucHom(\mathbb{D}_{\eps},F)$, where $\mathbb{R}_{\tau}\ucHom$ is the internal hom of $\Ho(\AffCtildetau)$, and $\mathbb{D}_{\eps} = \RuSpec(1[\eps])$, is the representable stack of the dual numbers $1[\eps]$ over $\cC$.\\

But the ratio $\Delta A / |\frac{d\langle A \rangle}{dt}|$ above suggests that we express $\Delta A$, or $\Der_A(Y,-) = \Delta A(B)$, in terms of $TA$, where $Y = \RuSpec B$, $B \in \Comm(\CE)$. One way to achieve this is to have a functor $\Der_A(Y,-) \rarr \Der_{TA}(Y,-)$. We do have such a functor; the projection $\mathbb{D}_{\eps} \rarr *$ induces a section $\sigma: A \rarr TA$ in $\AffCEtildetauE$ inducing $d\sigma: \Der_A(Y,-) \rarr \Der_{TA}(Y,-)$. Now recall that for $f: F \rarr G$ a morphism of stacks, $Y = \RuSpec B \rarr F$, one defines the simplicial set of relative derivations in this context as the standard homotopy fiber of the natural morphism $df$, which is defined by $\Der_{F/G}(Y,M) = \Map_{Y/\AffCEtildetauE/G}(Y[M],F)$, hence the ratio $\Delta A / |\frac{d\langle A \rangle}{dt}|$ is properly interpreted as $d\sigma$, and is functorially represented pointwise by $\Der_{A/TA}(Y,-)$ in a compact fashion. Thus the bifunctor $\Der_{A/TA}$ represents the ratio $\Delta A/|\frac{d\langle A \rangle}{dt}|$.\\

Finally the time energy uncertainty principle's left hand side is provided by $\Der_{V/TV}(-,-) \times \Der_E(-,-)$ for $V \in \AffCEtildetauE$, or pointwise by $\Der_{V/TV}(Y,-) \times \Der_E(X,-)$, with $X  = \RuSpec A$, $A \in \Comm(\cC)$, $Y = \RuSpec B$, $B \in \Comm(\CE)$. In a non-quantum algebra treatment of the uncertainty relation, we want this product of functors to be non-trivial, and this can be stated as saying that the associated geometric realization functor is non-contractible, where for $F$ and $G$ two functors into $\SetD$, one defines the geometric realization functor $|F \times G|$ by:
\beq
|F \times G|(A,B) = |FA \times_{\SetD} GB| \cong |FA| \times |GB| \nonumber
\eeq
as a product of $CW$ complexes, and we ask that it be non-contractible.

\section{Branes, Entropy and energy}
\subsection{Energy and entropy}
Recall that we have defined branes to be factorizations $F = V \circ E$ of stacks $F \in \AffCtildetau$. Energy is that part of stacks $F$ that materializes concepts into the natural realm, and time is a manifestation of their dynamics, thus branes provide a modelization of natural phenomena starting from information. Focusing on energy for the moment, we limit our attention to information theoretic processes.\\

We posit that since energy provides a manifestation of information in the natural realm, it also has a signature that can be measured, something we colloquially also refer to as energy, but which is merely an attribute of $E$ as a functor. We denote by $\eps$ this measurable attribute of $E$ and if needed one denotes this by $\eps_E$. Various collections of functors we refer to as a system, each functor $F$ corresponding to a state of such a system, with accompanying energy $E$ in $F = V \circ E$ having an associated measurable energy $\eps$ that we refer to as the energy of such a state.\\

To a system of states, generally speaking, one can associate Shannon's entropy (\cite{S1}, \cite{S2}), which measures its information content. Since we work with energy states, such an entropy can be made a functional of energy (\cite{J1}, \cite{J2}, \cite{BN},\cite{L}), in the following sense. In those references it is shown that the statistical entropy, which is defined in terms of energy, is a special case of Shannon's entropy, and thus carries informational content. It is not statistical energy however we work with, for the simple reason that it uses probabilities, which in themselves are double truncations, probabilities being idealizations of statistics, and statistics themselves resulting from truncating $\infty$-categorical constructs, which are complete and which we work with. Moreover using probabilities puts one in the position of observer, who is lacking some information, whereas we work from the perspective of an immediate algebraic manifestation $\cC$ of the ground state of natural phenomena (\cite{RG5}), so that we can consider that at this point all is known and there is no need to work with probabilities. Recall that this is equivalent to placing ourselves from the perspective of a Source, that has access to a universe $\cK$ of all knowledge, within which $IRV$s are constructed. Nevertheless, the statistical entropy's formula can still be used in this context whereby the probabilities are now interpreted as contributions of each state of the overall system. Those can be expressed in terms of the energy of such states. To summarize our position, we work with a system resulting from an amalgam of different states, each with a given energy. This system has  a measurable entropy, illustrative of its underlying informational content. To fix ideas, we use $S = - \sum p \log_2 p$ with $p = e^{-\beta \eps}/Z$ and $Z = \sum e^{-\beta \eps}$, the sum being over all states, $\eps$ being energy and $\beta$ denoting a constant. To say that the entropy measures the informational content of our system is an illustration of a broader, abstract statement, namely that energy transmutes information into matter, a fact we have exploited by defining energy as a stack in $\AffCtildetau$.\\

\subsection{Fuzzy entropy}
Shannon's entropy is defined for precise values of the energy. In practice however energy is never really well-known and what we are dealing with instead is fuzzy values $\eps \pm \Delta \eps$, $\Delta \eps$ a measurable component of $\Delta E$.\\

When we write an expression such as $\eps \pm \Delta \eps$, it will be understood that we work with a truncated value of $\Delta \eps$, which can easily be implemented by considering uncertainties of the form $\Delta \eps = x \eps$ with $x \in \mathring{I}$. Those then are the energies we work with. Consider an ideal situation whereby we have a large but finite number of states $N$. We let:
\beq
p_{N,i,x} = \frac{e^{-\beta(\eps_i \pm \Delta \eps_i)}}{Z} \nonumber
\eeq
where $\Delta \eps = x \eps$, $x \in (0,1)$. This defines the proportional contribution of the $i$th state with a focus indexed by $x$, to a system of $N$ states. In order to have $\sum p = 1$, which in this situation reads $\sum_i \int_{-\eps_i}^{\eps_i} p_{N,i,x} = 1$, we therefore need, with $d_i = x \eps_i$:
\beq
Z = \sum_i \int_{-\eps_i}^{\eps_i} e^{-\beta(\eps_i + d_i)}  \nonumber
\eeq
Armed with these notations, we finally have for all $N$:
\beq
S_N = - \sum_i \int_{-\eps_i}^{\eps_i} p_{N,i,x} \log_2 p_{N,i,x} \nonumber
\eeq

\section{Categories of motifs - Dharmakaya categories}
Of independent interest we discuss motifs (\cite{RG6}) and their various manifestations. This is an important topic since when we consider reconstruction functors in $\Fun(\cG, \Catinf)$ from $IRV$s (\cite{RG5}), there may be various objects being produced that share some common attributes, hence pointing to the existence of an underlying motif, itself being produced from $IRV$s, thereby presenting motifs as intermediate objects in the manifestation of natural phenomena. From that perspective we posit between $\cG$ and $\Catinf$ there are model categories of motifs that provide a signature for various objects of $\Catinf$.\\

Presently we focus with concepts with a same signature. Note that in a category one may have different classes of signatures, but we will focus on one generic signature that is fixed once and for all for illustrative purposes. Objects with the same signature behave alike, they exhibit the same patterns, their manifestations are the same in some regard. If $\Sigma$ denotes the class of those objects with a same signature $\sigma$, one has a class $\cO_{\Sigma}$ of objects that cannot distinguish them so if one has a morphism $f:X \rarr Y$ in $\Sigma$, and $W \in \cO_{\Sigma}$, $W$ not distinguishing between $X$ and $Y$ can be cast in the formalism of perceptions by askingthat the induced map of homotopy function complexes $\map(Y,W) \rarr \map(X,W)$ is a weak equivalence in $\SetD$. We will again refer to the class of such morphisms $f$ by $\Sigma$ for simplicity. Recall from \cite{Hi} that $\SetD$ is a left proper cellular model category, hence so is the category of simplicial presheaves over $\AffC$. The model categories of prestacks and stacks are both left Bousfield localizations originating from $\sPr(\AffC)$, hence by \cite{Hi} Proposition 4.1.1 are left proper cellular as well. Thus $\AffCtildetau$ is a left proper cellular model category, and again by Proposition 4.1.1 one can take a left Bousfield localization of this model category with respect to a class $\Sigma$ of morphisms therein. We take $\Sigma$ to be a class of morphisms between objects with a same signature as introduced above. It is worth recalling the standard definitions for a left Bousfield localization: if $f:A \rarr B$ is in $\Sigma$, a morphism between objects that share the same signature, an object $W$ is called $\Sigma$-local if it is fibrant, and the induced map $f^*: \map(B,W) \rarr \map(A,W)$ is a weak equivalence in $\SetD$. The class $\cO_{\Sigma}$ mentioned above therefore corresponds to the class of $\Sigma$-local objects. Those objects $W$ cannot detect objects $A$ and $B$ with a same signature, and this is what is signified by saying objects of that nature have the same manifestations from the perspective of other objects. From there one defines $\Sigma$-local equivalences to be those maps $f:X \rarr Y$ such that $\map(Y,W) \rarr \map(X,W)$ is a weak equivalence for any $\Sigma$-local object $W$. Thus what we have obtained is that $L_{\Sigma}\AffCtildetau$ is a model category structure on the underlying category of $\AffCtildetau$ for which weak equivalences are $\Sigma$-local equivalences. Upon taking the homotopy category of this category we get the category of motifs for objects of $\Sigma$.\\

\bigskip
\footnotesize
\noindent
\textit{e-mail address}: \texttt{rg.mathematics@gmail.com}.

\end{document}